%% file: main.tex
\documentclass[11pt,twoside]{amsart}
\input{math_preamble}

\usepackage{todonotes}

\usepackage[backend=biber,style=alphabetic]{biblatex}
 
\title{Entire curves in surface bundles over tori}
\author{Anthony Nguyen}
\AtEndDocument{{\footnotesize
		\textsc{School of Mathematics, University of Minnesota, Minneapolis, MN, US} \par
		\textit{E-mail address}: \texttt{nguy5343@umn.edu} \par		
}}

\begin{document}
\begin{abstract}	
		Cattalani recently initiated the study of Ahlfors currents in symplectic topology. We combine these with the symplectic pairs of Bande-Kotschick to show that some non-trivial surface bundles over tori have pseudoholomorphic entire curves. In particular the Kodaira-Thurston manifold has an entire curve for every $J$ tamed by an invariant form.  
\end{abstract}
	\maketitle

	\section{Introduction}
	A nonconstant holomorphic map $\C\to X$ is called an \textbf{entire curve} and a complex manifold $(X,J)$ is {\bf Brody hyperbolic} if it has no entire curves. This is equivalent to Kobayashi hyperbolicity for compact complex manifolds by Brody \cite{Brody78}.
	Both notions are extensively studied in complex geometry \cite{Langbook,Kobayashibook,Siu_survey,diverio_survey,Demailly_survey}. 
	
	These definitions of hyperbolicity extend to non-integrable structures: an almost complex manifold $(X,J)$ is called
	{\bf hyperbolic} if every $J$-holomorphic map
	\(\mathbb C\to X\) is constant. Otherwise, $(X,J)$ is {\bf non-hyperbolic} if it contains a $J$-holomorphic entire curve. One can furthermore consider whether a fixed symplectic form $\omega$ imposes either condition for all $\omega$-tame almost complex structure. This leads to the following definition: 
	\begin{definition}\label{def:hyperbolic}
		Let $(X,\omega)$ be a symplectic manifold.
		We say $(X,\omega)$ is {\bf hyperbolic} if, for every
		\(\omega\)-tame almost complex structure \(J\), all \(J\)-holomorphic
		maps $\C\to X$
		are constant.
		We say $(X,\omega)$ is {\bf non-hyperbolic} if, for every
		$\omega$-tame almost complex structure \(J\), there exists a
		nonconstant \(J\)-holomorphic map
		$\C\to X$. 
	\end{definition}
	The first result in this direction is due to Bangert \cite{Bangert1998Existence}. In contrast to the fact that tori never contain pseudoholomorphic rational curves, he proved the following:
	\begin{theorem}[Bangert \cite{Bangert1998Existence}]
		\label{thm:Bangert}
		There exists an entire $J$-holomorphic curve $\C\to T^{2n}$ for every almost complex structure $J$ tamed by the standard symplectic form $\w_\text{std}$.
	\end{theorem}
	
	Following Bangert's result, Gromov \cite{Spaces_and_Questions2000} proposed development of a theory of symplectic (non-)hyperbolicity. 
	His motivations for such an endeavor concern the study of symplectic structures on $(T^4,\w)$; his overarching goal is to show that every symplectic structure $\w$ on $T^4$ is isomorphic to the standard symplectic structure, as is known for $\CP^2$ from the study of rational curves\footnote{ We also see a similar story for rational and ruled symplectic 4-manifolds \cite{mcduff1990structure, LiLiu1995Ruled, LalondeMcDuff1996Classification, taubes1996swgr}.}.
	If his program is viable, its machinery should surely lend valuable insights into the symplectic structures of broad classes of (aspherical) 4-manifolds, insights that rational curves miss.
	
	Major steps toward Gromov's vision were undertaken by Cattalani \cite{spencerhyperbolic}, who initiated the study of Ahlfors currents in the symplectic category. With these objects he is able to  generalize Bangert's theorem by proving symplectic non-hyperbolicity for certain manifolds which are global products. We review his theorem in \autoref{section:prelim_results} along with all relevant definitions. 
	
	\subsection{Main theorem}
	We establish symplectic non-hyperbolicity of certain \textit{non-trivial} bundles by utilizing Ahlfors currents and the theory of symplectic pairs, introduced by Bande-Kotschick \cite{BandeKotschick2006}. 

\begin{theorem}\label{thm:main}
	Let $M$ be a closed, connected 4-manifold. Suppose that a symplectic pair $(\w_F,\w_B)$ is compatible with a positive genus surface bundle 
	$$\Sigma_g \to M \xrightarrow{\pi} T^2, \qquad g \geq 1.$$
	Then for the symplectic form\footnote{The notation implies $\w_F$ is a symplectic form on the fibers $\Sigma_g$ and $\w_B$ is the pullback of an area form on the base.} 
	$\w=\w_F + \w_B$, the space $(M,\w)$ is non-hyperbolic.
\end{theorem}
We refer the reader to \autoref{def:symplectic_pair} and \autoref{def:compatible_pair} in \autoref{section:prelim_results} for the definition of a compatible symplectic pair. The general argument of the theorem is similar to that of Cattalani. We use the symplectic pair to construct expanding Lagrangians, and we show the existence of holomorphic disks with prescribed Lagrangian boundary. Then we show that these disks satisfy the required estimates to approximate Ahlfors currents.

\subsection{Application to $T^2$-bundles over $T^2$} Geiges \cite{geiges1992symplectic} studied the existence of symplectic structures within each case of the Sakamoto-Fukuhara and Ue  \cite{SakamotoFukuhara1983,Ue1990} geometric classification of $T^2$-bundles over $T^2$. See \autoref{append:t^2_over_t^2} for an overview.
	
As an application of \autoref{thm:main}, we demonstrate symplectic non-hyperbolicity of invariant forms for new classes in \autoref{append:t^2_over_t^2}. In particular, 
\autoref{cor:KT-nonhyperbolicity} establishes case (b) and \autoref{cor:sol(3)-nonhyperbolicity} completes cases (e),(f),(g),(h). 
Combined with Bangert's seminal result for cases (a) and (c)\footnote{Bangert only explicitly proves case (a), but case (c) is then automatic since one can project an entire curve from a finite torus cover.}, this handles almost all classes of $T^2$-bundles over $T^2$. The only remaining case (d) is discussed in \autoref{conj:case_d}.

\begin{corollary}\label{cor:KT-nonhyperbolicity}
	Let $J$ be an almost complex structure on the Kodaira-Thurston manifold $K$ tamed by a left-invariant symplectic form. Then there exists an entire $J$-holomorphic curve $\C \to K$. 
\end{corollary}

\autoref{cor:KT-nonhyperbolicity} addresses the first part of \cite[Question 1]{spencerhyperbolic}. In light of \autoref{thm:main}, it follows by showing that every left-invariant symplectic form on $K$ decomposes as a sum of a compatible symplectic pair. Similarly, we have:

\begin{corollary}\label{cor:sol(3)-nonhyperbolicity}
	Let $M$ be a $T^2$-bundle over $T^2$ with geometric type $\Nil^3\times \R$ or $\Sol^3\times \R$. Let $J$ be an almost complex structure which is tamed by a left-invariant symplectic form on $M$. Then there exists an entire $J$-holomorphic curve $\C \to M$.
\end{corollary}

	\subsection*{Organization of the paper}
	In \autoref{section:prelim_results} we review relevant background on symplectic pairs and Ahlfors currents. In \autoref{section: Constructions} we construct the setting of the proofs; we use symplectic pairs to build Lagrangians and standard arguments to obtain pseudoholomorphic disks. 
	We then prove \autoref{thm:main} in \autoref{section:ahlfors} using Ahlfors currents. Finally in \autoref{section:T2bundles} we prove \autoref{cor:KT-nonhyperbolicity} and \autoref{cor:sol(3)-nonhyperbolicity}. 
	In \autoref{append:t^2_over_t^2} we recall the classification of $T^2$-bundles over $T^2$. 
	
	\subsection*{Acknowledgements}
	I am very grateful to my advisor Tian-Jun Li for his support and guidance. This work was supported by a GAANN Fellowship.

	\section{Preliminaries}
	\label{section:prelim_results}
	
	\subsection{Symplectic Pairs}
	It is natural to consider relaxing the non-degeneracy condition on a symplectic form. But only one degenerate form is often not enough to give non-trivial geometric structures. However, we notice from Donaldson's work \cite{donaldson2006twoformsfourmanifoldsellipticequations} that a pair or triple of 2-forms instead may yield interesting geometry.
	
	Symplectic pairs were introduced by Bande-Kotschick \cite{BandeKotschick2006}.
	Analogous to the story for symplectic forms, versions of Darboux theorem, Moser stability \cite{Bande2004Stability}, neighborhood theorems \cite{Her2015Neighborhood}, and contractibility of the space of compatible almost complex structures \cite{Her2018_ACS} also hold for symplectic pairs. We only consider them in four dimensions.

	\begin{definition} 
		\label{def:symplectic_pair}
		A \textbf{symplectic pair} on a 4-manifold $M$ is a pair of nontrivial closed 2-forms $(\w_1,\w_2)$ of constant and complementary ranks such that: 
		\begin{enumerate}[label=(\roman*)]
			\item $\w_1$ is a symplectic form on every leaf of $\ker \w_2$.
			\item $\w_2$ is a symplectic form on every leaf of $\ker \w_1$.
		\end{enumerate}
		The pair induces a symplectic form on $M$ given by 
	$\w^+ = \w_1 + \w_2$ and $\w^- = \w_1 - \w_2$.
	\end{definition}
	Manifolds with symplectic pairs are always symplectic, but they carry much stronger topological restrictions than general symplectic manifolds. For example, if a closed manifold $M^4$ admits a pair $(\w_1,\w_2)$, then $(M,\w^+) = (M,\w_1+\w_2)$ is minimal \cite{BandeGhiggini2020}. Also if $M$ admits a symplectic pair then it has a symplectic structure for both choices of orientation.

	\subsection{Compatibility with a bundle structure}	
	\begin{definition}
		\label{def:compatible_pair}
		A symplectic pair $(\w_F,\w_B)$ on $M$ is \textbf{compatible} with a surface bundle
		\[
		\Sigma_g \to M \xrightarrow{\pi} \Sigma_h
		\]
		if the kernel foliation of $\w_B$ is the fiber foliation of the surface bundle, i.e. 
		\[
		\ker \w_B = \ker d\pi 
		\]
	\end{definition}
	
	\begin{remark}
	Note that the existence of a symplectic pair on a total space does not automatically imply the existence of a \textit{compatible} symplectic pair for every bundle structure on that space. When compatible pairs do exist, they yield a flat connection and parallel transport preserves the symplectic form on the fibers.
	\end{remark}
	In general, a surface bundle over a surface
	$
	\Sigma_g \to M \xrightarrow{\pi} \Sigma_h
	$
	has a compatible symplectic pair if and only if $\pi$ admits a flat area-preserving connection \cite[Section 3.1]{BandeKotschick2006}. However, when $g\geq 3$, any $
	\Sigma_g \to M \xrightarrow{\pi} \Sigma_h
	$
	can be made flat by fiber summing with a trivial bundle \cite{Bowden_2011,Kotschick_Morita2005}. That is,
	 $
	 M\#_{\Sigma_g}(\Sigma_{g}\times \Sigma_{k})  \to \Sigma_{h+k}
	 $
	 can be made to support a compatible symplectic pair.
	
	For sphere fibers, Bande-Ghiggini \cite[Theorem 1.4]{BandeGhiggini2020} give a complete classification. They proved that 
	$
	S^2 \to M \to \Sigma_h
	$
	has a compatible symplectic pair for both the trivial and non-trivial bundles when $h\geq 1$. When $h=0$, they proved that only the trivial bundle has a compatible symplectic pair.  
	 
	In the case $g=h=1$, every oriented $T^2$-bundle over $T^2$ carries a symplectic pair \cite{BandeKotschick2006} and furthermore, examples of compatible symplectic pairs also exists. \cite{kahn2004symplectictorusbundlesgroup} gives a characteristic class criteria for the existence of a fiberwise symplectic form coming from restriction of a symplectic form $\Omega$ on $M$. Such a form can then be used to construct a compatible symplectic pair, but the symplectic pair does not necessarily sum to $\Omega$. A natural question is then:
	\begin{center}
		\textit{ When does a prescribed symplectic form on a total space necessarily decompose as the sum of a compatible symplectic pair?}
	\end{center}

\subsection{Hyperbolicity}	
As mentioned in the introduction, Cattalani \cite{spencerhyperbolic} initiated the study of Ahlfors currents in symplectic topology.
\begin{definition}
	Let $(X,J)$ be an almost complex manifold. A 2-current 
	$T:\Omega^2(X) \to \R$
	is an \textbf{Ahlfors current} if there exists a sequence of $J$-holomorphic disks $D_i\subset X$ such that
	\[
	\lim_{i\to \infty} \frac{D_i}{\area(D_i)} = T\quad\qand\quad \lim_{i\to \infty} \frac{\len(\partial D_i)}{\area(D_i)} = 0.
	\] 
\end{definition}
These objects are particularly useful because their existence is equivalent to the existence of entire pseudoholomorphic curves. The classical Ahlfors lemma \cite{ahlfors1935theorie} implies every entire curve yields an Ahlfors current, and Duval \cite{Duval2008SurLeLemme} established the opposite direction. Furthermore, such currents are much more tractable than directly dealing with entire curves.
\begin{proposition}(\cite{ahlfors1935theorie}, \cite{Duval2008SurLeLemme})
	Let $(X,J)$ be a compact almost complex manifold. Then $(X,J)$ contains an Ahlfors current if and only if it contains a non-constant $J$-holomorphic map $u: \C \to X$. 
\end{proposition}

Alongside the study of the moduli space and bubbling behavior of Ahlfors currents, Cattalani uses them to show the following existence result in certain global products:
\begin{theorem}[ \cite{spencerhyperbolic} Theorem 1.4]
	Let $(X,\w)$ be a closed connected symplectically aspherical manifold so that 
	\begin{enumerate}[label=(\alph*)]
		\item $X$ retracts onto a closed Lagrangian submanifold $L_X$ and
		\item the homological Dehn function $\delta_X$ of $\pi_1(X)$ grows at most quadratically 
	\end{enumerate}
	Then, there exists a nonconstant $J$-holomorphic map $u:\C \to X\times T^2$ for every almost complex structure tamed by $\w\oplus \w_{std}$. 
\end{theorem}

	\begin{remark}
		Fuchs-Yang \cite{fuchs2026holomorphicplanesmtimesmathbbt2} have observed that the hypothesis of Cattalani's theorem can be relaxed to just requiring the existence of a weakly exact Lagrangian. Similar arguments in this paper can also recover this version of Cattalani's theorem.
	\end{remark}

	\begin{remark}Outside of Ahlfors currents, existence results such as \cite{LiWu2012Note} and
		\cite{NNP-N26} are modifications of Bangert's arguments instead. The former used almost-K\"ahler cuts and the latter obtained a curvature criteria for symplectic non-hyperbolicity as a corollary of an almost-K\"ahler splitting theorem.

	\end{remark}

	\section{Lagrangian tori and pseudoholomorphic disks}
	\label{section: Constructions}
	
	In this section we construct the setting needed for the proof of \autoref{thm:main}. We first use symplectic pairs to construct the necessary Lagrangian tori and then we prove existence of the pseudoholomorphic disks that will eventually be used to estimate Ahlfors currents.
    
	Let $(M,\w)$ be as in the assumptions of \autoref{thm:main}. Then \[\Sigma_g \to M \xrightarrow{\pi} T^2, \qquad g \geq 1\]
	where $\w$ is from a compatible symplectic pair $(\w_F,\w_B)$. 
	Denote the kernel foliations by
	\begin{align*}
	\V := \ker \w_B \qand \H := \ker \w_F. 
	\end{align*}
	For each $R$, take a finite covering $T^2_R \to T^2$ where
	\[
	T_R^2 = \R^2 / (n_R\Z)^2, \quad n_R >2R 
	\]
	and consider the pullback bundle 
	$\pi_R:M_R \to T_R^2$.
	Note that $M_R$ is aspherical. 
	Denote again by $\w$ the pullback symplectic form on $M_R$. 

    \subsection{The Lagrangian tori $L_R$ with $\pi_2(M_R,L_R)\simeq \Z$}
	
	Fix an essential embedded loop $\gamma \subset \Sigma_g$. 
	The map $\pi_1(\Sigma_g) \mono\pi_1(M_R)$ is injective because $\pi_2(T^2)=0$.
	Choose a disk in  the universal cover $\R^2$ which embeds into $T_R^2$.
	Denote the image of the disk in $T_R^2$ by $B_R$. 
	Horizontal parallel transport yields a product trivialization 
	$\pi^{-1}_R(B_R(0)) \cong \Sigma_g \times B_R(0)$. 
	Set 
	\[S_R:= \partial B_R(0) \qand L_R:= \gamma \times S_R \subset M_R\]
	and choose $x_0\in \gamma$. The horizontal disk 
	\[D_{R,x_0}:= \{x_0\} \times B_R(0)\]
	has boundary on $L_R$ so we denote its relative homotopy class by 
	$\beta_R \in \pi_2(M_R,L_R)$. 
	
	\begin{lemma}
		\label{lem:finite_cover_constructions}
		The pair $(L_R,\beta_R)$ has the following properties. 
		\begin{enumerate}[label=(\roman*)]
			\item $L_R$ is an embedded Lagrangian torus in $(M_R,\w)$. 
			
			\item There are constants $0<c_0\leq C_0$ independent of $R$ such that
			\[c_0R^2 \leq \w(\beta_R) \leq C_0R^2.\] 
			
			\item Any disk bounding $L_R$ is generated by $\beta_R$, that is, \[\pi_2(M_R,L_R) \cong \Z \langle \beta_R\rangle \]
		\end{enumerate}
	\end{lemma}
	\begin{proof}
		At each point of $L_R$, the tangent space is spanned by a vertical vector tangent to $\gamma$ and a horizontal vector tangent to $S_R$. We know $\w_F$ vanishes on horizontal vectors and $\w_B$ vanishes on vertical vectors, hence $\w|_{L_R} = \w_F|_{L_R} + \w_B|_{L_R} = 0$.   
		Since the relative class is horizontal, 
		\[\int_{\beta_R}\w_F = 0 ,\qquad \int_{\beta_R} \w = \int_{\beta_R}\w_B = \int _{B_R(0)}\eta \]
		where $\eta$ is the area form on $T^2$. 
		Then (ii) follows from the area-growth estimate for disks with the Euclidean area. 
		Now we show (iii). 	
		Since $\pi_2(M_R)=0$, from the long exact sequence for relative homotopy we know that
		$\pi_2(M_R,L_R) \simeq \ker \{ \pi_1(L_R) \xrightarrow{i_*} \pi_1(M_R)\}$. Consider an arbitrary element 
		$
		a[\gamma] + b[S_R] \in \pi_1(L_R).
		$
		Since $\gamma$ is essential and $\Sigma_g$ has no torsion, the image of $[\gamma]$ in $\pi_1(M_R)$ has infinite order and cannot map to zero, as $\pi_1(\Sigma_g)\to \pi_1(M_R)$ is injective. 
		By construction we have $i_*[S_R]=0$, which implies 
		\[
		\ker \{ \pi_1(L_R) \to \pi_1(M_R)\} \simeq \Z \langle [S_R]\rangle 
		\]
		This implies that 
		$\pi_2(M_R,L_R) = \Z \langle\beta_R\rangle$ since $\partial \beta_R = [S_R]$.
	\end{proof}
	
	\subsection{Existence of pseudoholomorphic disks}
	Now we consider the existence of pseudoholomorphic disks for the almost complex structures on the finite covers $M_R$. We adapt the argument in \cite[Lemma 2.26]{spencerhyperbolic} to a local setting.

	Choose a complex structure $j_R$ on $T_R^2$ and a $\w_F$-compatible complex structure on $\V = \ker \w_B $. Consider the horizontal lift of $j_R$ to $TM_R$. Let be $J_{0,R}$ the almost complex structure defined from the vertical complex structure and the horizontal lift of $j_R$. Then $J_{0,R}$ is $\w$-compatible, preserves the splitting $TM_R=\V\oplus \H$, and is a product almost complex structure on $\Sigma_g\times B_R(0)$. 
	For any almost complex structure $J$ on $M_R$, let $\M_{0,1}(J,[D])$ denote the space of $J$-holomorphic disks homologous to $D\in \pi_2(M_R,L_R)$ with one marked point on the boundary. Let 
	\[
	\ev_{J}\colon \M_{0,1}(J,[D]) \to L
	\]
	be the evaluation map taking a disk to its marked point. 
	
	\begin{lemma}\label{lem:deg_model_disks}
		Every $J_{0,R}$-holomorphic disk in the class $\beta_R$ is Fredholm regular and the evaluation map 
		$
		\ev_{J_{0,R}}:\M_{0,1}(J_{0,R},\beta_R) \to L_R
		$
		has \[\deg(\ev_{J_{0,R}}) =  1\mod 2.\]
	\end{lemma}
	\begin{proof}
		Let $u$ be a $J_{0,R}$-holomorphic disk in the class $\beta_R$. 
		As $J_{0,R}$ preserves the splitting, we have
		\[
		u^*\w_F(\partial_s, \partial_t) = \w_F(v,J_{0,R}|_\V (v))\geq 0\qquad (s,t)\in \D\]
		for any vertical vector $v \in \V$ along the image of $u$. Note that 
		\[
		\int_D u^* \w_F = \w_F(\beta_R) = 0
		\]
		because a representative disk is horizontal. This implies that $v=0$ hence the image of $u$ lies in one horizontal leaf.
		In the trivialization over $B_R(0)$, each horizontal leaf is $\{x\}\times B_R(0)$. The boundary of $u$ lies in $\gamma \times S_R$, so the constant point in the fiber satisfies $x\in \gamma$. Then $u(\D,\partial \D)$ is one of the horizontal disks 
		\[
		D_{R,x} = \{x\}\times B_R(0), \qquad x\in \gamma 
		\]
		up to reparametrization. Along a horizontal disk, the linearization of the Cauchy-Riemann operator is the direct sum $D_u=D_{u_1}\oplus D_{u_2}$ for
		\begin{enumerate}[label=(\roman*)]
			\item $D_{u_1}$ the standard linearized operator at a constant map $u_1:(\D,\partial \D) \to \{x\}$.
			
			\item $D_{u_2}$ the standard linearized operator for
			$u_2:(\D,\partial \D) \to \{T^2_R,S_R\}$.
		\end{enumerate}
		Both operators are surjective, hence every disk in the class $\beta_R$ is regular. After quotienting by reparametrization, the disk $u$ is determined by $x\in \gamma$ and the boundary marked point is given by its image in $S_R$. Thus
		$
		\M_{0,1}(J_{0,R},\beta_R) \cong \gamma \times S_R = L_R.
		$
		This shows that $\deg(\ev_{J_{0,R}}) = 1\mod 2$ as the evaluation map is the identity under this identification. 
	\end{proof}

	\begin{proposition}
		\label{prop:J-disk_prescribed_beta}
		Let $J_R$ be any almost complex structure on $M_R$ tamed by $\w$. Then there is a nonconstant $J_R$-holomorphic disk 
		\[u_R: (\D,\partial \D) \to (M_R,L_R)\]
		representing $\beta_R$. Every such disk has fixed area 
		\[
		\area(u_R) = \int_\D u_R^*\w = \w(\beta_R).
		\]
	\end{proposition}
	\begin{proof}
		We show disks in the class $\beta_R$ are simple.
		In \cite{lazzarini2000existence} and \cite[Theorem A]{ Lazzarini2011Relative} Lazzarini shows that for any disk $u$ in the class $\beta_R$ there exists finitely many simple disks $\{v_i\}_{i=1}^N$ and $\{m_i \}_{i=1}^N\in \N$ with the same Lagrangian boundary such that
		\[
		u(\D) = \bigcup_{i=1}^N v_i(\D) \quad\qand\quad  [u] = \sum_{i=1}^N m_i[v_i] \quad \text{ in } H_2 (M_R,L_R;\Z)
		\]
		Since $[v_i]\in \pi_2(M_R,L_R)$, we know from \autoref{lem:finite_cover_constructions} that $[v_i] = k_i\beta_R \in \pi_2(M_R,L_R)$.
		Note that the homological decomposition yields
		\begin{align*}
		\w(\beta_R) = \area(u)  &= \sum_{i=1}^N m_i \area(v_i) = \sum_{i=1}^N m_ik_i \w(\beta_R) 
		\implies
		1 = \sum_{i=1}^N m_ik_i
		\end{align*}
		which implies that $N=1$ and $m_1=k_1=1$. Hence we obtain a simple disk $v$ such that $[v]=\beta_R$ and $\area(u)=\area(v)$ and $u(\D)=v(\D)$. This proves that disks in $\beta_R$ are simple. 
		Disk bubbling cannot occur because every nonconstant disk component with boundary on $L$ has area at least $\w(\beta_R)$, so there is exactly one component. Sphere bubbling does not occur because $\pi_2(M_R)=0$. 
		Thus the moduli spaces $\M_{0,1}(J,\beta_R)$ are compact for a fixed tamed almost complex structure and for a compact one-parameter family.
		
		Next, we first consider a special almost complex structure. Let $J_{0,R}$ be the $\w$-compatible almost complex structure defined above which is a product almost complex structure on $\Sigma_g\times B_{R}(0)$. 
		The evaluation map 
		\[
		\ev_{J_{0,R}}\colon \M_{0,1}(J_{0,R},\beta_R) \to L_R  
		\]
		has degree 1 mod 2 by \autoref{lem:deg_model_disks}. 
		Consider a generic tame $J_{1,R}$ and choose a generic path of tame almost complex structures from $J_{0,R}$ to $J_{1,R}$. The parameterized moduli space is a compact cobordism between the moduli spaces 
		$\M_{0,1}(J_{0,R},\beta_R)$ and $\M_{0,1}(J_{1,R},\beta_R)$. 
		This provides a bordism between $\ev_{J_{0,R}}$ and $\ev_{J_{1,R}}$. Hence 
		\[\deg(\ev_{J_{1,R}}) = 1 \mod 2\]
		So the moduli space for $J_{1,R}$ is nonempty. For an arbitrary $\w$-tamed $J_R$ a Gromov compactness argument with approximating generic tame structures yields the desired result. 
		Then for such a $J_R$-holomorphic disk $u_R$ representing $\beta_R$, the condition  $\area(u_R) = \w(\beta_R)$ is immediate.
	\end{proof}
	
	\section{Entire curves from Ahlfors currents}
	\label{section:ahlfors}
	
	In this section, we prove \autoref{thm:main}. 
	
	For an arbitrary $\w$-tamed $J$ on $M$, lift to $J_R$ on $M_R$. 
	Choose a $J_R$-holomorphic disk into $M_R$ from \autoref{prop:J-disk_prescribed_beta}. Let $\tilde u_R\colon(\D,\partial \D) \to (\tilde M, \tilde L_R)$ denote the lift of such $u_R$ to the universal cover, which is a bundle 
	$\tilde \pi \colon \tilde M \to \R^2$ and
	$\tilde L_R=\tilde \gamma \times S_R$
	for $\tilde \gamma$ a lift of $\gamma$. 
	We show that given a sequence of holomorphic disks in the universal cover $\tilde M$, we can extract an Ahlfors current on the original space $M$.
	
	Let $(X,g)$ be a Riemannian manifold and $d$ the induced metric. Let $S\subset X$ be nonempty. Denote the distance function to $S$
	\[
	d_S: X \to [0,\infty), \qquad d_S(p) := \inf_{x\in S}\{d(p,x)\}
	\]
	For $a\in (0,\infty]$ denote 
	\[
	N_a(S) := d_S^{-1}\left( [0,a) \right)
	\]
	the set of points of distance at most $a$ from $S$. 
	We use the following general version of Cattalani's disk shrinking lemma. 
	\begin{lemma}[c.f. \cite{spencerhyperbolic} Lemma 2.22]
		\label{lem:disk_shrinking}
		Let $(\D,g)$ be a disk with a Riemannian metric $g$. There exists an embedded disk $\D'\subset \D$ with rectifiable boundary such that $\partial \D'\subset N_a(\partial \D)$ and 
		\[
		\frac{\len(\partial \D')}{\area(\D')} \leq \frac{1}{a} \frac{\area(N_a(\partial \D))}{ \area(\D-N_a(\partial \D))}
		\]
	\end{lemma}
	\begin{proof}
		Fix a Riemannian disk $(\D,h)$ and consider the rescaled metric $\hat h = \frac{1}{a^2} h$ where $0<a<\infty$. Notice that
		\[
		\len_{\hat h} = a^{-1}\len_{h}
		\:
		\implies 
		\:
		N_1^{\hat h}(\partial D) = N_a^{h}(\partial D)
		\qand \area_{\hat h} = a^{-2} \area_h. 
		\]
		Then we apply \cite[Lemma 2.22]{spencerhyperbolic}  to $(\D, \hat h)$ to obtain the inequality
		\[
		a\frac{\len_{ h}(\partial \D')}{\area_h(\D')} 
		=
		\frac{\len_{\hat h}(\partial \D')}{\area_{\hat h} (\D')} 
		\leq  
		\frac{\area_{\hat h}(N_1(\partial \D))}{ \area_{\hat h}(\D-N_1(\partial \D))} 
		= 
		\frac{a^{-2}}{a^{-2}} \frac{\area_ h(N_a(\partial \D))}{ \area_ h(\D-N_a(\partial \D))} 
		\]
		which implies for $(\D,h)$ that
		\[
		\frac{\len_h(\partial \D')}{\area_h (\D')} \leq \frac{1}{a} \frac{\area_ h(N_a(\partial \D))}{ \area_ h(\D-N_a(\partial \D))} 
		\]
	\end{proof}
	By \cite[Remark 2.23] {spencerhyperbolic} this lemma also applies to metrics induced by nonconstant pseudoholomorphic maps, e.g. $h_R=u^*_Rg_J$ considered below. 
	
	\begin{proposition}
		\label{lem:ahlfors_current_existence}
		Assume the preceding setup.
		For $R\to \infty$, let $u_R\colon (\overline \D,\partial \overline \D) \to (M_R,L_R)$
		be any holomorphic disk representing $\beta_R$, and let \[\tilde u_R:(\overline \D,\partial \overline \D) \to \tilde M\] be their lifts to the universal cover. 
		Then there is a constant $C$ and smoothly embedded subdisks $\D'_R\subset  \interior \overline \D$ such that 
		\[
		\frac{\len (\tilde u_R |_{\partial \D'_R})}{\area (\tilde u_R|_{\D'_R})} \leq \frac{C}{R} \to 0.
		\]
		After projection, a normalized subsequence converges to a nonzero Ahlfors current on $M$. 
	\end{proposition}
	
	\begin{proof} 
		The main steps are as follows.  
		\begin{enumerate} [label=(\arabic*)]
			\item  We first show that there are subdisks whose lower area bound grows quadratically.
			\item We show this lower area bound and \autoref{lem:finite_cover_constructions} implies that the ratio 	
			\[
			\frac{\area(N_R(\partial \D))}{ \area(\D-N_R(\partial \D))}
			\]
			is uniformly bounded, independent of $R$. 
			\item We apply \autoref{lem:disk_shrinking}, to get pseudoholomorphic disks with length/area ratio of order $O(R^{-1})$.
			
			\item We then project to $M$ and normalize to obtain an Ahlfors current.
		\end{enumerate}

		\textbf{(1)} For any tamed $J$ denote the induced metric $g_J$ on $M$ and denote $g_0$ a flat metric on $T^2$. The induced tangent map of the bundle projection $\pi\colon (M,g_J) \to (T^2,g_0)$ has bounded operator norm since
		\[
		\| d\pi_p \| := \sup_{v\in T_pM, v\neq 0}\frac{|d\pi_p(v)|}{|v|_{g_J}}, \qquad \qquad \|d\pi\| := \sup_{p\in M} \| d\pi_p \| < \infty
		\]
		and $M$ is compact. On each finite cover $\pi_R\colon M_R\to T_R^2$, we use the lifted metrics of $g_J$. The supremum, denoted $K$, bounds $\|d\pi_R\|$ for all $R$ since $M_R \to M$ is a local isometry. 
		Recall that $B_R$ is the embedded disk in $T^2_R$ with $S_R=\partial B_R$ from \autoref{section: Constructions}.
		Choose a smooth bump function $\chi_R:T^2_R \to [0,1]$ with $\supp(\chi_R) \subset B_{R/2}$ and
		\[
		\chi_R = 1 \qquad \text{ on }B_{R/4}.
		\] 
		Let $\mu$ denote an area form on $T^2_R$ and set
		\[
		\mu_R = \chi_R\mu \qand  \alpha_R = \pi_R^* \mu_R.
		\]
		$\alpha_R$ is closed because $d\alpha_R = d\pi_R^*(\mu_R)=\pi_R^*(d\mu_R)=0$ and the support $\supp(\alpha_R)\subset M_R$ is disjoint from $L_R=\gamma \times S_R$, hence $\alpha_R|_{L_R}=0$. 
		Since 
		$u_R$ and $D_{R,x_0}$ represent the same relative homotopy class $\beta_R$, we choose such a homotopy
		\[
		H:[0,1]\times \overline \D \to M_R, \qquad H([0,1]\times \partial \D)\subset L_R.
		\]
		Stokes theorem gives
		\[
		\int_\D u^*_R\alpha_R - \int_{D_{R,x_0}}\alpha_R = \int_{[0,1]\times \D}H^*d\alpha_R - \int_{[0,1]\times \partial \D} H^*\alpha_R = 0
		\]
		Then it follows that
		\[
		\int_\D u^*_R \alpha_R = \int_{D_{R,x_0}} \alpha_R = \int_{B_R} \mu_R \geq \area(B_{R/4}) = \frac{\pi R^2}{16}
		\]
		Now let 
		$
		V_R = (\pi_R \circ u_R)^{-1}(B_{R/2})\subset \D
		$
		and note that $u^*_R \alpha_R$ vanishes outside $V_R$. We have 
		\[
		|u^*_R\alpha_R(v,w)|  = |u_R^*\pi_R^* \mu_R(v,w)| \leq K^2 |v\wedge w| 
		\]
		so that 
		\[
		\frac{\pi R^2}{16} \leq \int_{V_R}u^*_R\alpha_R \leq  \int_{V_R} |u^*_R\alpha_R| \leq K^2 \area(u_R|_{V_R})
		\]
		hence for a constant $c$ we have
		$\area(u_R|_{V_R}) \geq cR^2.$ 
		
		\textbf{(2)} Let $h_R = \tilde u_R^* \tilde g_J$ be the pullback metric and $d_{h_R}$ the induced distance. For $\rho_R = R/4K$ now consider $N_{R}:=N_{\rho_R}(\partial \D)$ a neighborhood of the domain disk boundary. 
		$V_R$ is disjoint from $N_R$. The preceding lower bound implies
		\[
		\area(u_R) - \area_{h_R}(N_R) =  \area_{h_R}(\overline \D\setminus N_R) \geq \area_{h_R}(V_R) \geq cR^2
		\]
		Using $\area(u_R)\leq C_0R^2$ from  \autoref{lem:finite_cover_constructions} we obtain
		\[
		\area(u_R) -  \area_{h_R}(N_R) \geq \delta \area (u_R), \qquad \delta = \frac{c}{C_0}
		\]
		For sufficiently large $C_0$, we have that $0<\delta<1$ and 
		\[
		\frac{\area_{h_R}(N_R)}{\area(u_R) - \area_{h_R}(N_R)} \leq \frac{1-\delta}{\delta}.
		\]
		
		\textbf{(3)} Now by \autoref{lem:disk_shrinking}
		we obtain subdisks $\D'_R\subset \D$ with rectifiable boundary such that 
		\[
		\frac{\len (\tilde u_R |_{\partial \D'_R})}{\area (\tilde u_R|_{\D'_R})} \leq \frac{C}{R} \to 0.
		\]
		
		Consider the projection of the disks to $M$.
		A subsequence of the normalized currents of integration converges to an Ahlfors current by \cite[Lemma 26.14]{Simon1983}.
		\end{proof}

\begin{proof}[Proof of \autoref{thm:main}]
		Fix an arbitrary $\w$-tamed almost complex structure $J$ on $(M,\w)$ and a sequence $\{R_k\} \to \infty$. For each $R_k$ we construct a finite cover $M_{R_k}$, Lagrangian torus $L_{R_k}\subset M_{R_k}$, and primitive class $\beta_{R_k}\in \pi_2(M_{R_k},L_{R_k})$ as in \autoref{lem:finite_cover_constructions}. 
		
		By \autoref{prop:J-disk_prescribed_beta} the lifted almost complex structure admits a pseudoholomorphic disk $(\D,\partial \D) \to (M_{R_k},L_{R_k})$ in the class $\beta_{R_k}$. 
		Lift such disks to the universal cover.

		Thus \autoref{lem:ahlfors_current_existence} guarantees an Ahlfors current on $(M,J)$, and Duval's quantitative Brody Lemma \cite{Duval2008SurLeLemme} gives a nonconstant $J$-holomorphic map $\C \to M$. 
\end{proof}

\section{Left-invariant forms on $T^2$-bundles over $T^2$}
	\label{section:T2bundles}
	
	In this section, we prove \autoref{cor:KT-nonhyperbolicity} and \autoref{cor:sol(3)-nonhyperbolicity}. 
	To prove \autoref{cor:KT-nonhyperbolicity} it suffices to show the following:
\begin{proposition}
		\label{prop: kt_pair_splitting}
		Every left-invariant symplectic form $\w$ on the Kodaira-Thurston manifold $K$ admits a decomposition into a symplectic pair 
		$
		\w = \w_F + \w_B
		$
		compatible with a bundle structure
		$
		T^2 \to K \xrightarrow{\pi} T^2.
		$
\end{proposition}
	
\begin{proof}[Proof]
	Let $G=\operatorname{Heis}_3(\mathbb R)\times\mathbb R$ so that $K=\Gamma\backslash G$ and $\w$ is any left-invariant symplectic form. Following Geiges \cite{geiges1992symplectic},
	consider fibrations $\pi_1,\pi_2:K\to T^2$ induced by
	\[
	(x,y,z,t)\mapsto(x,t) \qand
	(x,y,z,t)\mapsto(y,t)
	\]
	respectively. These correspond to the first and third descriptions of $K$ as a torus bundle in \cite[Section~5(b)]{geiges1992symplectic}. Let $F_1,F_2$ denote torus fibers of these fibrations.
	Geiges showed that for every $a\in H^2(K;\mathbb R)$ we have the following:
	\begin{enumerate}
		\item If both $a(F_1) = 0$ and $a(F_2) = 0$ then $a^2=0$.
		
		\item $a$ has a symplectic form representative if and only if $a^2\neq 0$.
	\end{enumerate}
	Now let $a=[\w ]$. Taken together, (1) and (2) imply both fibrations cannot simultaneously have $\w$ vanishing on its fibers. It follows that $\int_{F_i}\w\neq 0$ for at least one $i\in\{1,2\}$, so choose such an index and set $\pi=\pi_i$.
	Its vertical distribution $\mathcal V=\ker d\pi$ is invariant and
	$\w|_{\mathcal V}$ is nondegenerate.
	Consequently, the complement $\H=\mathcal V^\w$ is an invariant rank-two
	distribution, and there is a symplectic orthogonal splitting
	\[
	TK=\mathcal V \oplus \H.
	\]
	Choose a translation-invariant area form $\sigma$ on the base torus, and let $\beta = \pi^* \sigma.$ Notice $\pi$ is induced by a Lie group homomorphism, so the lift of $\beta$ to $G$ is left-invariant. 
	The form $\beta$ is closed,
	$\ker\beta=\mathcal V$, and $\beta|_{\H}$ is nondegenerate.
	Since $\w|_\H$ and $\beta|_\H$ are nonzero two-forms on a two-dimensional distribution, there is a unique smooth non-vanishing function $c:K\to\mathbb R$
	such that
	\[
	\w|_{\H}=c\,\beta|_{\H}.
	\]
	The lifts of $\w$, $\beta$, and $\H$ are all left-invariant,
	so the lift of $c$ is a left-invariant function on $G$ and is therefore
	constant. Thus $c$ is constant on $K$. 
	Define
	\[
	\w_B=c\beta,\qquad \w_F=\w-c\beta.
	\]
	Constancy of $c$ implies $d\w_B=d\w_F=0$.
	By construction,
	\[
	\w_F|_{\H}=0,\qquad
	\w_B|_{\mathcal V}=0,\qquad
	\w_F|_{\mathcal V}=\w|_{\mathcal V},\qquad
	\w_B|_{\H}=\w|_{\H}.
	\]
	Moreover, for $v\in \V$ and $h\in\H$,
	\[
	\w(v,h)=0,\qquad \beta(v,h)=0,
	\]
	so both $\w_F$ and $\w_B$ have zero mixed terms.
	Since their nonzero restrictions are nondegenerate, it follows
	that
	\[
	\ker\w_F=\H,\qquad \ker\w_B=\mathcal V.
	\]
	Both forms therefore have constant rank two. Their kernel distributions
	are integrable because the forms are closed, and each form restricts
	symplectically to the leaves of the kernel of the other. Thus $(\w_F,\w_B)$ forms a symplectic pair compatible with the bundle structure for $K$, and their sum is the prescribed
	form $\w$.
\end{proof}
We now restate and prove \autoref{cor:sol(3)-nonhyperbolicity}. 
\begin{corollary*}[\autoref{cor:sol(3)-nonhyperbolicity}]
	Let $M$ be a $T^2$-bundle over $T^2$ with geometric type $\Nil^3\times \R$ or $\Sol^3\times \R$. Let $J$ be an almost complex structure tamed by a left-invariant symplectic form on $M$. Then there exists an entire $J$-holomorphic curve $\C \to M$.
\end{corollary*}
	\begin{proof}
		Consider $T^2$-bundles $M$ of geometric type $\Sol^3\times \R$.
		These are the bundles cases (g) and (h) in \autoref{append:t^2_over_t^2}. 
		In light of \autoref{thm:main}, we just need to show every left invariant symplectic form on $M$ can be decomposed with a symplectic pair compatible with $\pi$.
		\cite[Theorem 18, Case (g,h)]{BandeKotschick2006} shows that symplectic pairs on $M$ exist. 
		Furthermore by \cite[Theorem 2]{geiges1992symplectic} any prescribed symplectic form on $M$ always has $\w$ non-vanishing on its fibers.
		Therefore the construction in \autoref{prop: kt_pair_splitting} can essentially be repeated verbatim.
		
		The same applies to bundles with type  $\Nil^3\times \R$ with $b_1=2$ \cite[Theorem 1]{geiges1992symplectic} and for the $b_1=3$ case this was \autoref{prop: kt_pair_splitting}. One can alternatively use the previous result \autoref{cor:KT-nonhyperbolicity} and note that cases (e) and (f) have a Kodaira-Thurston finite cover, so $J$-holomorphic entire curves descend via projection. The nonconstant property follows since the covering map is a local homeomorphism. 
	\end{proof}	
	
\begin{conjecture}[Case (d)]
	\label{conj:case_d}
	Let $M$ be a $T^2$-bundle over $T^2$ with geometric type $\Nil^4$. Let $J$ be an almost complex structure which is tamed by a left-invariant symplectic form on $M$. Then there exists a nonconstant entire $J$-holomorphic curve $\C \to M$.
\end{conjecture}
This is the only remaining case in the classification of $T^2$-bundles over $T^2$. 
In this setting, the torus fibration does not have symplectic torus fibers (see \autoref{append:t^2_over_t^2}). So although they admit pairs \cite[Theorem 18 Case (d)]{BandeKotschick2006}, they are never compatible with the fibration. Nonetheless, we still expect symplectic non-hyperbolicity. 
One heuristic reason might be \cite[Lemma 2.5, Conjecture 4.9]{Ho_Li_2012}, which shows the Lagrangian bundle structure is preserved by Luttinger surgery and predicts that all $T^2$-bundles over $T^2$ are related by such surgery. Since cases (a),(b), and (d) have Lagrangian bundle structure and are related under surgery, one might expect entire curves to be preserved if they avoid the Lagrangian fiber. 

In general it would be interesting  to study entire $J$-holomorphic curves with further qualitative control, such as those which either always avoid or pass through prescribed sets. 

\appendix
\section{Classification of orientable $T^2$-bundles over $T^2$}
\label{append:t^2_over_t^2}

The following is the classification of $T^2$-bundles over $T^2$ and their Thurston geometries, due to \cite{SakamotoFukuhara1983} and \cite{Ue1990} respectively. The matrices describe the monodromy corresponding to generators of $\pi_1(T^2)$, and pairs of integers $(m,n)$ represent the Euler class. 

Geiges \cite[Theorem 1]{geiges1992symplectic} shows that outside of cases (b) and (d) the fibration of $M$ as a torus bundle is unique.
\cite[Theorem 2]{geiges1992symplectic} shows that in cases (a) and (b) a given invariant symplectic form may or may not be compatible with a selected fibration. For case $(d)$ the symplectic structure is never compatible with the fibration. For all other cases the symplectic structure is always compatible with the fibration.

\begin{center}
	\begin{tabular}{|c|c|l|c|}
		\hline
		& $b_1$ & Monodromy \& Euler class & Geometric type \\
		\hline\hline
		
		(a) & 4
		& $\{I,I,(0,0)\}=T^4$
		& $\mathbb{R}^4$ \\
		\hline
		
		(b) & 3
		& $\{I,I,(m,n)\},\quad (m,n)\neq(0,0)$
		& $\Nil^3\times\mathbb{R}$ \\
		\hline
		
		(c) & 2
		& $\left\{
		\begin{pmatrix}0&-1\\1&-1\end{pmatrix},
		I,(0,0)
		\right\}$
		& $\mathbb{R}^4$ \\
		& &
		$\left\{
		\begin{pmatrix}0&-1\\1&-1\end{pmatrix},
		I,(-1,0)
		\right\}$
		& \\
		& &
		$\left\{
		\begin{pmatrix}0&-1\\1&0\end{pmatrix},
		I,(0,0)
		\right\}$
		& \\
		& &
		$\left\{
		\begin{pmatrix}0&-1\\1&0\end{pmatrix},
		I,(-1,0)
		\right\}$
		& \\[3pt]
		& &
		$\left\{
		\begin{pmatrix}1&-1\\1&0\end{pmatrix},
		I,(0,0)
		\right\}$
		& \\
		& & $\{-I,I,(0,0)\}$ & \\
		& & $\{-I,I,(-1,0)\}$ & \\
		\hline
		
		(d) & 2
		& $\left\{
		\begin{pmatrix}1&\lambda\\0&1\end{pmatrix},
		I,(m,n)
		\right\},\quad \lambda\neq0,\ n\neq0$
		& $\Nil^4$ \\
		\hline
		
		(e) & 2
		& $\left\{
		\begin{pmatrix}-1&\lambda\\0&-1\end{pmatrix},
		I,(m,n)
		\right\},\quad \lambda\neq0$
		& $\Nil^3\times\mathbb{R}$ \\
		\hline
		
		(f) & 2
		& $\left\{
		\begin{pmatrix}1&\lambda\\0&1\end{pmatrix},
		-I,(m,n)
		\right\},\quad \lambda\neq0$
		& $\Nil^3\times\mathbb{R}$ \\
		\hline
		
		(g) & 2
		& $\{C,I,(m,n)\},\quad
		|\tr C|\geq3,\ C\in\SL_2(\mathbb{Z})$
		& $\Sol^3\times\mathbb{R}$ \\
		\hline
		
		(h) & 2
		& $\{C,-I,(m,n)\},\quad
		\tr C\geq3,\ C\in\SL_2(\mathbb{Z})$
		& $\Sol^3\times\mathbb{R}$ \\
		\hline
	\end{tabular}
\end{center}

\printbibliography

\end{document}

%% file: math_preamble.tex
\usepackage{hyperref}
\hypersetup{
	colorlinks=true,
	allcolors=gray
}
\usepackage[margin=1.35in]{geometry}
\usepackage{slashed}
\usepackage{amssymb}
\usepackage{amsmath}
\usepackage{amsthm}
\usepackage{mathrsfs}
\usepackage{dutchcal}
\usepackage{color}
\usepackage{wrapfig}
\usepackage{import}
\usepackage{xifthen}
\usepackage{pdfpages}
\usepackage{transparent}
\usepackage{eucal}
\usepackage{tensor}
\usepackage{parskip}
\usepackage{mathtools}
\usepackage{extpfeil}
\usepackage{thmtools}{}
\usepackage{tikz}
\usepackage{tikz-cd}
\usetikzlibrary{arrows.meta}

\tikzset{
	symbol/.style={%
		draw=none,
		every to/.append style={%
			edge node={node [sloped, allow upside down, auto=false]{$#1$}}}
	}
}{}
\usepackage[bbgreekl]{mathbbol}
\usepackage[shortlabels]{enumitem}
\usepackage{setspace}
\DeclareMathOperator{\area}{area}
\DeclareMathOperator{\interior}{int}
\DeclareMathOperator{\len}{length}

\DeclareMathOperator{\Sol}{Sol}
\DeclareMathOperator{\Nil}{Nil}

\DeclareMathOperator{\tr}{tr}

\DeclareMathOperator{\SL}{SL}

\DeclareMathOperator{\ev}{ev}

\DeclareMathOperator{\supp}{supp}

\newcommand{\Z}{\mathbb{Z}}
\newcommand{\R}{\mathbb{R}}

\newcommand{\mono}{\hookrightarrow}

\renewcommand{\phi}{\varphi}
\newcommand{\w}{\omega}

\newcommand{\C}{\mathbb{C}}
\newcommand{\N}{\mathbb{N}}

\newcommand{\D}{\mathbb{D}}

\renewcommand{\H}{\mathcal{H}}

\newcommand{\M}{\mathcal{M}}
\newcommand{\V}{\mathcal{V}}

\newcommand{\CP}{\mathbb{CP}}

\renewcommand{\tilde}{\widetilde}

\renewcommand{\implies}{\Rightarrow}

\newcommand{\qand}{\quad\text{and}\quad}

\declaretheorem{theorem}
\declaretheorem[numbered=no,name=Theorem]{theorem*}
\declaretheorem[numbered=no,name=Corollary]{corollary*}
\declaretheorem[sibling=theorem]{lemma}
\declaretheorem[sibling=theorem]{proposition}
\declaretheorem[sibling=theorem]{corollary}
\declaretheorem[sibling=theorem]{conjecture}

\declaretheorem[style=remark,sibling=theorem]{remark}

\declaretheoremstyle[
spaceabove=6pt, spacebelow=6pt,
headfont=\normalfont\bfseries,
notefont=\mdseries, notebraces={(}{)},
bodyfont=\normalfont,
postheadspace=1em
]{definition}
\declaretheorem[style=definition,sibling=theorem]{definition}